# Monotone Domain Decomposition Iterative Method and its Convergence for a Nonlinear Integro-Differential Equation of Volterra Type

Myong-Gil RIM and Dong-Hyok KIM

*Faculty of Mathematics, **Kim Il Sung** university, Pyongyang, D.P.R. Korea*

***Abstract***: We apply the monotone domain decomposition iterative method to a nonlinear integro-differential equation of Volterra type and prove its convergence. To do this, by adding a term in both sides of the original equation we make a linear equation to get a monotone domain decomposition iterative scheme and prove the existence, uniqueness and convergence of iterative solutions.

***Keyword***: Domain decomposition, Monotone iterative scheme, convergence

***MSC***:     65R20, 65N55

## 1.  Introdution

The monotone domain decomposition iterative method is one of the important methods for solving partial differential equations. Domain decomposition methods are said to be a collection of techniques based on the principle of Divide and Conquer. These methods have been developed for solving partial differential equations on domains in 2D and 3D. These methods solve a boundary value problem by splitting it into smaller boundary value problems on subdomains and iterating to coordinate the solution between adjacent subdomains. The problems on the subdomains are independent, which makes domain decomposition methods suitable for parallel computing. In [1, 8, 9] is suggested the monotone domain decomposition iterative method for a type of nonlinear parabolic differential equation and proved the convergence. In overlapping domain decomposition methods, the subdomains overlap by more than the interface. Overlapping domain decomposition methods include the Schwarz alternating method and the additive Schwarz method. In [2, 5, 6] they proved the convergence of some multiplicative and additive Schwarz methods for inequalities which contain contraction operators.

However, as far as the authors know, the domain decomposition methods have not been used in integro-differential equations. In [3, 4, 7] are suggested the methods of subsolution and supersolution for a type of nonlinear integro- differential equation of Volterra type and proved their convergence.

In paper we suggest a monotone domain decomposition iterative method for the following nonlinear integro-differential equation of Volterra type

$$
\left.
\begin{aligned}
& u_t - Lu = f(t,\ x,\ u) + \int_0^t g_0(t,\ x,\ s,\ u(t,\ x),\ u(s,\ x))ds, \quad D \\
& Bu \equiv \alpha_0(t,\ x)\partial u/\partial \nu + \beta_0(t,\ x)u = h(t,\ x), \qquad\qquad S \\
& u(0,\ x) = u_0(x), \qquad\qquad\qquad\qquad\qquad\qquad \Omega
\end{aligned}
\right\}
$$

(1)

and prove its convergence.

In the remainder of the paper, we provide monotone domain decomposition iterative scheme for (1) and then prove its convergence.

## 2. Monotone domain decomposition iterative scheme for (1)

In (1) let $\Omega \subset \mathbf{R}^N$ be a bounded domain with a smooth boundary, $D = (0,\ T) \times \Omega$, $S = (0,\ T) \times \partial\Omega$ and $L$ a uniformly elliptic operator given by

$$
Lu \equiv \sum_{i,\ j=1}^n a_{ij}(t,\ x)\partial^2 u/\partial x_i \partial x_j + \sum_{j=1}^n b_j(t,\ x)\partial u/\partial x_j .
$$

Here, the matrix $(a_{ij})$ is positive definite in $\overline{D}$. We suppose that

$\alpha_0 \equiv \alpha_0(t,\ x)$, $\beta_0 \equiv \beta_0(t,\ x)$ are continuous on $S$, and $\alpha_0 \geq 0$, $\beta_0 \geq 0$, $\alpha_0 + \beta_0 > 0$.

**Definition:** Let $g_0(\,\cdot\,,\ \eta_1,\ \eta_2)$ be nondecreasing in $\eta_2$. Then a function $\tilde{u} \in C(\overline{D}) \bigcap C^{1,\,2}(D)$ is called a ***supersolution*** of equation (1) if it satisfies the following inequality.

$$
\left.
\begin{aligned}
& \tilde{u}_t - L\tilde{u} \geq f(t,\ x,\ \tilde{u}) + \int_0^t g_0(t,\ x,\ s,\ \tilde{u}(t,x),\ \tilde{u}(s,\ x))ds, \quad D \\
& \alpha_0 \partial\tilde{u}/\partial \nu + \beta_0\tilde{u} \geq h(t,\ x), \qquad\qquad\qquad\qquad\qquad S \\
& \tilde{u}(0,\ x) \geq u_0(x), \qquad\qquad\qquad\qquad\qquad\qquad\quad \Omega
\end{aligned}
\right\} \qquad (2)
$$

Similarly, $\hat{u}$ is called a ***subsolution*** if it satisfies the reversed inequality.

Suppose that there exists a pair of ordered supersolution and subsolution $\tilde{u}$, $\hat{u}$.

Suppose that $\Omega$ is composed of two subdomains, that is, $\Omega = \Omega_1 \bigcup \Omega_2$.

Let $J = \langle \hat{u}, \tilde{u} \rangle := \{ u \in C(\overline{D}) \mid \hat{u} \le u \le \tilde{u}, \overline{D} \}$.

We assume that $f$ is a Hölder continuous function in $D \times J$ with the exponent $\alpha$ and a $C^1$-function in $J$ and $g_0$ is a Hölder continuous function in $D \times [0, T] \times \boldsymbol{R}^2$ and $g_0(\cdot, \eta_1, \eta_2)$ is continuously differentiable when $\eta_1, \eta_2 \in J$.

Also suppose that $g_0$ satisfies a Lipschitz condition, that is, for each $r > 0$ there exists a constant $K_0 \equiv K_0(r) > 0$ such that

$$| g_0(\cdot, \eta_1, \eta_2) - g_0(\cdot, \eta_1', \eta_2') | \le K_0(|\eta_1 - \eta_1'| + |\eta_2 - \eta_2'|), \ |\eta_i| \le r, \ |\eta_i'| \le r, \ i = 1, 2 \quad (3)$$

Let $g(t, x, u(t, x)) \equiv \int_0^t g_0(t, x, s, u(t, x), u(s, x)) ds$ and for each $\forall t, s \in (0, T], \ x \in \Omega$ we define

$$\underline{b}_0(t, x, s) = \sup \left\{ -\frac{\partial g_0}{\partial \eta_1}(t, x, s, \eta_1, \eta_2), \ \eta_1, \eta_2 \in J \right\},$$

$$\underline{b}(t, x) = \int_0^t \underline{b}_0(t, x, s) ds$$

Without loss of generality we may assume that $\underline{b} \in C^\alpha(D)$ for each $u, \mathsf{v} \in J$ with $u \ge \mathsf{v}$ in $D$, the nondecreasing property of $g_0(\cdot, \eta_1, \eta_2)$ in $\eta_2$ and definity of $\underline{b}_0, \ \underline{b}$ imply that

$$g_0(\cdot, u(t, x), u(s, x)) - g_0(\cdot, \mathsf{v}(t, x), \mathsf{v}(s, x)) =$$
$$= [g_0(\cdot, u(t, x), u(s, x)) - g_0(\cdot, \mathsf{v}(t, x), u(s, x))] +$$
$$+ [g_0(\cdot, \mathsf{v}(t, x), u(s, x)) - g_0(\cdot, \mathsf{v}(t, x), \mathsf{v}(s, x))] \ge -\underline{b}_0[u(t, x) - \mathsf{v}(t, x)],$$

that is,

$$g(t, x, u(t, x)) - g(t, x, \mathsf{v}(t, x)) \ge -\underline{b}(t, x)[u(t, x) - \mathsf{v}(t, x)].$$

From the above results for each $u \in J$, $(t, x) \in D$ the function

$$F_1(t, x, u) = (\underline{c} + \underline{b})u + f(t, x, u) + g(t, x, u)$$

is nondecreasing in $u$. Here $\underline{c}(t, x) = \sup\{-f_u(t, x, u) ; \hat{u} \le u \le \tilde{u}\}$.

By letting $c = \underline{c} + \underline{b}$ and adding the term $cu$ on both sides of the first equation of (1), we have

$$\mathbf{L}_c[u] \equiv (\partial/\partial t - L + c)[u] = F_1(t, x, u), \ D.$$

Hence we make the following monotone domain decomposition iterative scheme.

$$u_{11}^{(0)} = u_{21}^{(0)} = \hat{u}, \quad u_{12}^{(0)} = u_{22}^{(0)} = \tilde{u}, \ \overline{D}$$

$$\left.\begin{array}{ll} \mathbf{L}_c[u_{1j}^{(n+1)}] = F_1(t,\ x,\ u_{1j}^{(n)}), & D_1 \\ Bu_{1j}^{(n+1)} = Bu_{2j}^{(n)}, & S_1 \\ u_{1j}^{(n+1)}(0,\ x) = u_0(x), & \Omega_1 \end{array}\right\},\ \ u_{1j}^{(n+1)} := u_{2j}^{(n)},\ \ \overline{D}\setminus\overline{D}_1$$

$$\left.\begin{array}{ll} \mathbf{L}_c[u_{2j}^{(n+1)}] = F_1(t,\ x,\ u_{2j}^{(n)}), & D_2 \\ Bu_{2j}^{(n+1)} = Bu_{1j}^{(n+1)}, & S_2 \\ u_{2j}^{(n+1)}(0,\ x) = u_0(x), & \Omega_2 \end{array}\right\},\ \ u_{2j}^{(n)} := u_{1j}^{(n)},\ \ \overline{D}\setminus\overline{D}_2$$

Then $u_{ij}^{(n)}$ is $n^{\text{th}}$ *iterative solution*.

## 3. Convergence of the above Iterative Scheme

**Lemma 1.**[2] *If* $\omega \in C(\overline{D}) \cap C^{1,2}(D)$ *satisfies*

$$\left.\begin{array}{ll} \omega_t - L\omega + c\omega \ge 0, & D \\ \alpha_0 \partial\omega/\partial\nu + \beta_0 \partial\omega \ge 0, & S \\ \omega(0,\ x) \ge 0, & \Omega \end{array}\right\}\ (\ \alpha_0 \ge 0,\ \beta_0 \ge 0,\ \alpha_0 + \beta_0 > 0,\ \ S\ )$$

*and* $c \equiv c(t,\ x)$ *is a bounded function in* $D$, *then we have*

$$\omega(t,\ x) \ge 0,\ \ \overline{D}.$$

**Lemma 2**[2] *Consider the linear parabolic boundary-value problem*

$$\left.\begin{array}{ll} u_t - Lu + cu = q(t,\ x), & D \\ \alpha_0(t,\ x)\partial u/\partial\nu + \beta_0(t,\ x)u = h(t,\ x), & S \\ u(0,\ x) = u_0(x), & \Omega \end{array}\right\} \qquad (4)$$

*The following two statements holds:*

① *Let* $\alpha_0 = 0$ *and* $q$ *be locally Hölder continuous locally in* $x$ *and uniformly in* $t$. *Then for any continuous* $h$, $u_0$ *which satisfy the compatibility condition*

$\beta_0(0,\ x)u_0(x) = h(0,\ x),\ \partial\Omega$ *the first boundary-value problem (4) has a unique solution* $u(t,\ x)$.

$$u(t,\ x) = J^{(1)}(t,\ x) + \int_0^t d\tau \int_\Omega G(t,\ x;\ \tau,\ \xi)q(\tau,\ \xi)d\xi + \int_0^t d\tau \int_{\partial\Omega} \frac{\partial\Gamma}{\partial\nu_\xi}(t,\ x;\ \tau,\ \xi)\psi(\tau,\ \xi)d\xi$$

②*Let* $\alpha_0 = 1$ *and* $q$ *be Hölder continuous in* $x$ *and Hölder continuous uniformly in* $\overline{D}$. *Then for any continuous function* $h$, $u_0$ *the secondary boundary-value problem (4) has a unique solution* $u(t,\ x)$ *which is Hölder continuous in* $x$.

$$u(t, \ x) = J^{(0)}(t, \ x) + \int_0^t d\tau \int_\Omega \Gamma(t, \ x; \tau, \ \xi) q(\tau, \ \xi) d\xi + \int_0^t d\tau \int_{\partial\Omega} \Gamma(t, \ x; \tau, \ \xi) \psi(\tau, \ \xi) d\xi$$

Here $J^{(1)}$, $J^{(0)}$, $\Gamma$, $G$, $\psi$ is givend in work[2].

**Theorem**  Let $f(\cdot, \eta)$, $g_0(\cdot, \eta_1, \eta_2)$ be $C^1$-functions in $\eta$ and $\eta_1$, respectively. Assume that $g_0(\cdot, \eta_1, \eta_2)$ is nondecreasing in $\eta_2 \in J$ and satisfies the Lipschitz condition (3). Then the sequences $\{u_{ij}^{(n)}\}$ given by the iterative scheme converges monotonically to the unique solution $u$ of (1).

**Proof.**  First, prove that the $n$ th iterative solution sequence $\{u_{ij}^{(n)}\}$ is well defined.

By Lemma 2 $q_{ij}^{(n)}(t, \ x) \equiv F_1(t, \ x, \ u_{ij}^{(n)}(t, \ x))$ is Hölder continuous[2] and hence $n$ th iterative solution sequence $\{u_{ij}^{(n)}\}$ is well defined.

Second, prove that $\{u_{ij}^{(n)}\}$ posses the monotone property, that is,

$$\hat{u} \le u_{11}^{(n)} \le u_{21}^{(n)} \le u_{11}^{(n+1)} \le u_{12}^{(n+1)} \le u_{22}^{(n)} \le u_{12}^{(n)} \le \tilde{u}, \ \ \overline{D} \tag{5}$$

For $n = 0$, $\hat{u} = u_{11}^{(0)} = u_{21}^{(0)}$, $\tilde{u} = u_{22}^{(0)} = u_{12}^{(0)}$.

Prove that $u_{11}^{(1)} \ge u_{21}^{(0)}$, $\overline{D}$.

$$(\partial_t - L + c) u_{11}^{(1)} = F_1(t, \ x, \ u_{11}^{(0)}) = F_1(t, \ x, \ u_{21}^{(0)}), \ \ D_1$$

$$(\partial_t - L + c) u_{21}^{(0)} \le F_1(t, \ x, \ u_{21}^{(0)}), \ D_1$$

$$(u_{11}^{(1)} - u_{21}^{(0)})_t - L(u_{11}^{(1)} - u_{21}^{(0)}) + c(u_{11}^{(1)} - u_{21}^{(0)}) \ge 0, \ \ D_1$$

$$B(u_{11}^{(1)} - u_{21}^{(0)}) = Bu_{11}^{(1)} - Bu_{21}^{(0)} = 0, \ \ \ S_1$$

$$u_{11}^{(1)}(0, \ x) - u_{21}^{(0)}(0, \ x) = u_0(x) - u_{21}^{(0)}(0, \ x) \ge 0, \ \ \ \Omega_1$$

By Lemma 1   $u_{11}^{(1)} \geq u_{21}^{(0)}$,  $\overline{D}_1$ .

Also because  $u_{11}^{(1)} = u_{21}^{(0)}$  in  $\overline{D} \setminus \overline{D}_1$ ,  $u_{11}^{(1)} \geq u_{21}^{(0)}$,  $\overline{D}$ .

In the same way  $u_{12}^{(1)} \leq u_{22}^{(0)}$, $u_{11}^{(1)} \leq u_{12}^{(1)}$,  $\overline{D}$ . Thus for  $n = 0$  equation (5) holds.

We Assume that the equation (5) holds for  $n - 1$ .
By the work[2]

$$\hat{u} \leq u_{11}^{(n)} \leq u_{11}^{(n+1)} \leq u_{12}^{(n+1)} \leq u_{12}^{(n)} \leq \tilde{u}, \ \ \overline{D}_1, \ \ \hat{u} \leq u_{21}^{(n)} \leq u_{21}^{(n+1)} \leq u_{22}^{(n+1)} \leq u_{22}^{(n)} \leq \tilde{u}, \ \ \overline{D}_2 .$$

Using  $u_{1j}^{(n+1)} = u_{2j}^{(n)}$  in  $\overline{D} \setminus \overline{D}_1$ , then  $\hat{u} \leq u_{11}^{(n)} \leq u_{11}^{(n+1)} \leq u_{12}^{(n+1)} \leq u_{12}^{(n)} \leq \tilde{u}, \ \ \overline{D}$ .

Prove that  $u_{11}^{(n+1)} \geq u_{21}^{(n)}$,  $\overline{D}$ .

$$u_{21}^{(n)} = u_{11}^{(n)} \leq u_{11}^{(n+1)} \ \text{in} \ D_1 \setminus D_2 \ \text{and} \ u_{11}^{(n+1)} = u_{21}^{(n)} \ \text{in} \ D_2 \setminus D_1 .$$

$$(\partial_t - L + c)(u_{11}^{(n+1)} - u_{21}^{(n)}) = F_1(t, \ x, \ u_{11}^{(n)}) - F_1(t, \ x, \ u_{11}^{(n-1)}) \geq 0 \ \text{in} \ D_{12} \equiv D_1 \bigcap D_2 .$$

$$Bu_{11}^{(n+1)} = Bu_{21}^{(n)} \ \text{in} \ S_1 \bigcap D_2 \ \text{and} \ Bu_{11}^{(n+1)} = Bu_{21}^{(n+1)} \geq Bu_{21}^{(n)} \ \text{in} \ S_2 \bigcap D_1 .$$

Hence  $Bu_{11}^{(n+1)} \geq Bu_{21}^{(n)}$,  $D \bigcap (S_1 \bigcup S_2)$  and  $u_{11}^{(n+1)}(0, \ x) - u_{21}^{(n)}(0, \ x) = 0$,  $\Omega_1 \bigcap \Omega_2$.

By Lemma 1  $u_{11}^{(n+1)} \geq u_{21}^{(n)}$,  $\overline{D}_{12}$ . Hence  $u_{11}^{(n+1)} \geq u_{21}^{(n)}$,  $\overline{D}$ .

In the similar way  $u_{21}^{(n)} \geq u_{11}^{(n)}$,  $u_{22}^{(n)} \geq u_{12}^{(n+1)}$,  $u_{12}^{(n)} \geq u_{22}^{(n)}$,  $\overline{D}$  and the equation (5) holds.

From the equation (5)   the limit  $\lim_{n \to \infty} u_{ij}^{(n)} = u_{ij}$,  $i, j = 1, \ 2$  exists and

$$u_{1j} = u_{2j} \ (j = 1, 2), \ \ u_{i1} \leq u_{i2} \ (i = 1, 2) .$$

By the work[2]  $u_1 = u_{i1}, \ u_2 = u_{i2}, \ i = 1, \ 2$  is a solution of nonlinear integro-parabolic equation of Volterra type (1).

The uniqueness of the solution is proved by the work[2].

Thus the theorem is proved.